\documentclass[conference]{IEEEtran}
\IEEEoverridecommandlockouts
\usepackage{cite}
\usepackage{lmodern}
\usepackage{amsmath}
\usepackage{graphicx}
\usepackage[colorlinks=true, allcolors=black]{hyperref}
\usepackage{amssymb}
\usepackage{amsfonts}
\usepackage{amsthm}
\usepackage{colortbl}
\usepackage{dsfont}
\usepackage{bm}
\usepackage{algorithm}
\usepackage{algorithmic}
\usepackage{dsfont}
\usepackage[table, dvipsnames]{xcolor}
\usepackage{tabularx}
\usepackage{lipsum}
\usepackage{mathtools}
\usepackage{listings}

\newcommand{\espE}{\mathbb{E} }
\newcommand{\Vdom}{ U }

\newtheorem{theorem}{Theorem}

\newtheorem{proposition}[theorem]{Proposition}

\newtheorem{definition}[theorem]{Definition}

\def\BibTeX{{\rm B\kern-.05em{\sc i\kern-.025em b}\kern-.08em
    T\kern-.1667em\lower.7ex\hbox{E}\kern-.125emX}}
\begin{document}

\title{A Gradient Method for Risk Averse Control of a PDE-SDE Interconnected System}

\author{\IEEEauthorblockN{Gabriel Velho}
\IEEEauthorblockA{\textit{\small Université Paris-Saclay,} \\
\textit{\small CentraleSupélec, CNRS} \\
\textit{\small Laboratoire des Signaux} \\ 
\textit{\small et Systèmes,} \\
\small Gif-sur-Yvette, France \\
\footnotesize gabriel.velho@centralesupelec.fr}
\and
\IEEEauthorblockN{Jean Auriol}
\IEEEauthorblockA{\textit{\small Université Paris-Saclay,} \\
\textit{\small CentraleSupélec, CNRS} \\
\textit{\small Laboratoire des Signaux} \\ 
\textit{\small et Systèmes,} \\
\small Gif-sur-Yvette, France \\
\footnotesize jean.auriol@centralesupelec.fr}
\and
\IEEEauthorblockN{Riccardo Bonalli}
\IEEEauthorblockA{\textit{\small Université Paris-Saclay,} \\
\textit{\small CentraleSupélec, CNRS} \\
\textit{\small Laboratoire des Signaux} \\ 
\textit{\small et Systèmes,} \\
\small Gif-sur-Yvette, France \\
\footnotesize riccardo.bonalli@centralesupelec.fr}
}

\maketitle

\begin{abstract}

In this paper, we design a risk-averse controller for an interconnected system composed of a linear Stochastic Differential Equation (SDE) actuated through a linear parabolic heat equation.
These dynamics arise in various applications, such as coupled heat transfer systems and chemical reaction processes that are subject to disturbances.
While existing optimal control methods for these systems focus on minimizing average performance, this risk-neutral perspective may allow rare but highly undesirable system behaviors. To account for such events, we instead minimize the cost within a coherent risk measure. Our approach reformulates the coupled dynamics as a stochastic PDE, approximates it by a finite-dimensional SDE system, and applies a gradient-based method to compute a risk-averse feedback controller. Numerical simulations show that the proposed controller substantially reduces the tail of the cost distribution, improving reliability with only a minor reduction in average performance.
\end{abstract}

\begin{IEEEkeywords}
Stochastic systems; Distributed parameter systems; Optimal control
\end{IEEEkeywords}


\section{Introduction}

Interconnected systems composed of parabolic PDEs coupled with lower--dimensional dynamics arise in many engineering applications, such as heat networks, chemical reactors, and fluid–structure interactions \cite{tang_state_2011, antonio_susto_control_2010, mohammadi_optimal_2015}. In particular, PDE–ODE systems coupled through their boundaries have been extensively studied, with stabilization methods ranging from finite-dimensional approximations \cite{auriol_late-lumping_2019} to backstepping designs \cite{vazquez_backstepping_2026}. However, in realistic settings, the finite-dimensional component is often subject to stochastic disturbances due to fluctuating environmental conditions or measurement uncertainty. Modeling these effects naturally leads to PDE–SDE interconnected systems, which capture both spatial diffusion and random perturbations \cite{touzi_optimal_2013, lew_risk-averse_2024}. A representative example is district heating\cite{talebi_review_2016}, where the PDE describes heat propagation in a pipeline and the SDE models the internal temperature dynamics of a building influenced by unpredictable external factors.

A first study on the stabilization and optimal control of such PDE–SDE systems was carried out in \cite{velho_optimal_2025}, but the objective there was to minimize the expected performance cost. While expectation--based (risk--neutral) control can yield good performance on average, it does not prevent rare but severe deviations, which correspond to system failures or excessive energy usage \cite{wang_risk-averse_2022}. In safety-critical settings, these rare events are precisely those that must be controlled.

Risk-averse stochastic control addresses this issue by minimizing a risk measure rather than a mean cost, explicitly reducing the likelihood of catastrophic outcomes \cite{shapiro_lectures_2009, chapman_risk-sensitive_2022, bonalli_first-order_2023}.
Effective computational methods exist for finite-dimensional SDEs, including gradient-based optimization of CVaR--type objectives \cite{miller_optimal_2017, velho_gradient_2025}, but no analogous framework currently exists for infinite-dimensional or PDE–SDE systems.

The work presented in this paper takes a step toward filling that gap, building directly on two of our previous contributions. These are: 1. the optimal control framework for parabolic PDE–SDE interconnected systems introduced in \cite{velho_optimal_2025}, which provides the foundation finite-dimensional approximation of interconnected systems, and 2. the gradient-based risk-averse control method developed in \cite{velho_gradient_2025}, which offers an effective strategy for optimizing coherent risk measures. By combining these two ideas, we obtain a tractable approach for risk-averse control of PDE–SDE systems: we reformulate the dynamics as a SPDE, approximate it with a finite-dimensional SDE, and apply the risk-averse gradient method to the reduced model. This procedure yields an approximate solution to the original infinite-dimensional risk-averse control problem. Establishing theoretical bounds on this approximation remains an interesting direction for future work. 

To our knowledge, this provides this yields the first methodology enabling risk-aware control of distributed systems under stochastic disturbances. Although built upon existing techniques, the full procedure requires nontrivial adaptations to handle the infinite-dimensional structure and the coupling, and goes beyond a direct plug-and-play combination of prior results.


The remaining of the paper is organized as follows: After recalling some standard notations, the problem under consideration is stated in Section \ref{PF}.  In Section \ref{CC}, we show how the coupled PDE–SDE system can be rewritten as a SPDE and then approximated by a finite--dimensional SDE. Section \ref{CS} presents the gradient--based method used to compute the risk-averse control law for the approximated system. Numerical results are provided in Section \ref{NR}.


\section{Problem formulation}\label{PF}

For a filtered probability space $(\Omega, \mathcal{F} \triangleq (\mathcal{F}_t)_{t \in [0,\infty)}, \mathbb{P})$, 
we assume that stochastic perturbations are due to a one--dimensional, $\mathcal{F}$--adapted Wiener process $W_t$. 

Let $T > 0$ be some given time horizon. For a Hilbert space $H$, we denote by 
$L_\mathcal{F}^2(0,T ; H)$ the set of square integrable processes $Z: [0,T]\times\Omega \to H$ that are $\mathcal{F}$--progressively measurable.

We denote $H^1(0,1) \triangleq \bigl \{u \in L^2(0,1 ; \mathbb{R}), \ \text{such that} \  u' \in L^2(0,1 ; \mathbb{R})  \bigr \}$. It is a Hilbert space with the scalar product
$$
< u , v >_{H^1} \triangleq \int_0^1 u(x) v(x) dx + \int_0^1 u'(x) v'(x) dx.
$$

\subsection{Control system}

We consider interconnected PDE--SDE systems of the form
\begin{equation}\label{eq:systeme_original_coupled}
\left\{ \begin{array}{l} \partial_t u(t,x) = \Delta u(t,x)  +  c u(t,x) \\ 
dX_t = \bigl( A X_t + B u(t,0) + r(t) \bigr) dt \\
\hspace{3em} + \bigl( C X_t + D u(t,0) + \sigma(t) \bigr) dW_t \\
\partial_x u(t,1) + \beta_1 u(t,1) = V(t) \\
\partial_x u(t,0) - \beta_0 u(t,0) =  M X_t  \\
u(0,x) = u_0(x), \quad X(0) = X_0
\end{array}  \right.
\end{equation}
in the time-space domain $[0, T] \times [0,1]$, where $T>0$. The state of the system is $(X_t, u(t,\cdot)) \in \mathbb{R}^d \times H^1(0,1)$. The operator $\Delta$ represents the Laplacian. The coefficient $c$ is a constant modeling the rate of growth ($c>0$) or decay ($c<0$) of the solution $u$ at each point in space. Coefficients $\beta_0$ and $\beta_1$ are positive, this ensures that the system is well posed \cite{brezis_functional_2011}. 
The matrices $A \in \mathbb{R}^{d\times d}, B \in \mathbb{R}^{d\times 1}, C \in \mathbb{R}^{d\times d}, D \in \mathbb{R}^{d\times 1}, M \in \mathbb{R}^{1 \times d}$ are also constant. The drifts $r(\cdot), \sigma(\cdot) \in L^\infty(0,T;\mathbb{R}^d)$ are known. The control input $V(t)$ takes values in $\mathbb{R}$. 
The considered domain of the operator $\Delta$ is
\begin{align*}
\text{Dom}(\Delta) \triangleq \bigl \{   u \in H^2(0,1), \ & \partial_x u(0) - \beta_0 u(0) = 0 \\
& \partial_x u(1) + \beta_1 u(1) = 0   \bigr \}.
\end{align*}
The initial data $(X_0, u_0)$ belongs to $\mathbb{R}^d \times \text{Dom}(\Delta)$ and is taken to be deterministic. We work with the state space $H^1(0,1)$ for the PDE component to ensure continuity at the boundary, which is guaranteed by the trace theorem for $H^1(0,1)$. Note that positive values of $c$ and positive eigenvalues of $A$ may introduce growth or instability in the respective subsystems.

\begin{figure}[ht!]
    \centering
    \includegraphics[width=0.9\linewidth]{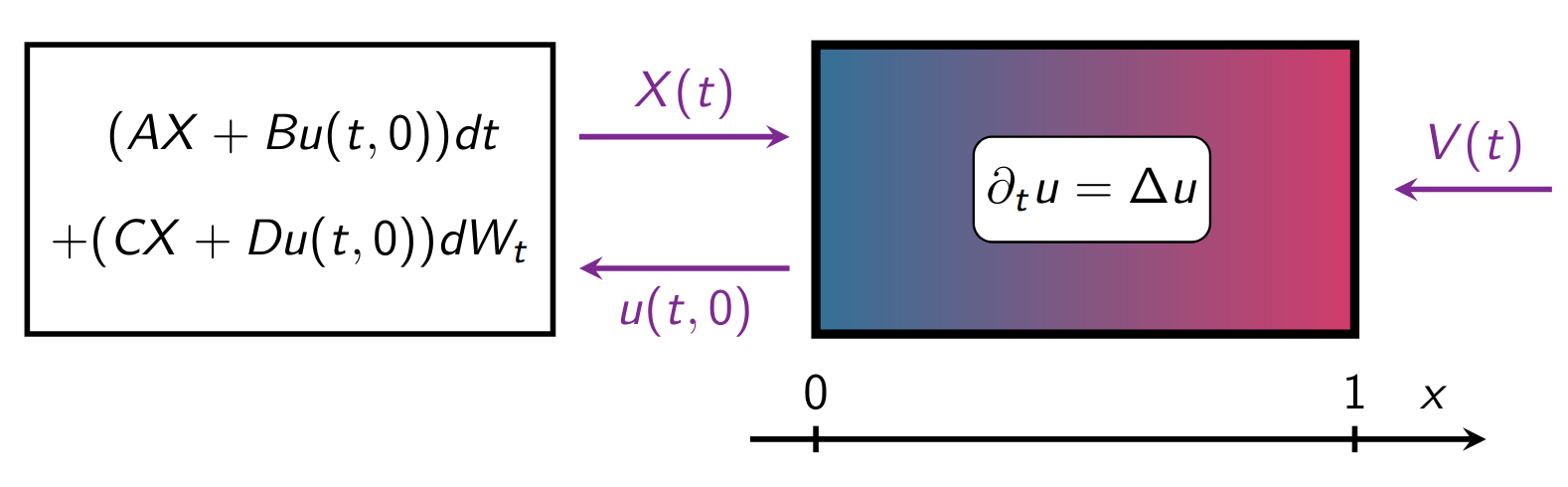}
    \caption{Structure of an interconnected parabolic PDE -- SDE}
    \label{fig:fig_parabolic_SDE_interconnected}
\end{figure}

Systems of the form \eqref{eq:systeme_original_coupled} naturally arise when a diffusive medium such as a heat conductor interacts with a lumped stochastic subsystem. A typical example is the temperature regulation of a building through a pipe transmitting heat, where external factors like fluctuating outdoor temperatures or varying sunlight exposure introduce stochastic disturbances \cite{talebi_review_2016}. In these settings, the Robin--type boundary conditions model a heat flux between the systems proportional to the temperature difference. 
The proof of the well--posedness of system~\eqref{eq:systeme_original_coupled} can be found in \cite[Prop. 2]{velho_optimal_2025}.

\subsection{Cost and optimization problem}

Our goal is to design a control law $V(\cdot)$ that ensures reliable regulation of the system. To this end, we consider the quadratic cost functional
\begin{equation}\label{eq:original_LQ_cost_constant}
J(V) \triangleq \int_0^T X_t^\top Q X_t  dt + X_T^\top G X_T + r \Vert u(\cdot) \Vert_{\mathcal{U}}^2 ,
\end{equation}
where $X_t$ evolves according to \eqref{eq:systeme_original_coupled}, and $\mathcal{U}$ is the control space. The weighting matrices $Q$ and $G$ belong to $\mathcal{S}_q^+$, and $r>0$ is a given scalar. The precise definition of $\mathcal{U}$ is presented in subsection \ref{CHOICE_CONTROL}, once the necessary additional context has been introduced. Note that $J(V)$ is a random variable, due to the stochastic nature of the dynamics.


In \cite{velho_optimal_2025}, we proposed a method for minimizing the \textit{expected value} of $J(V)$. While this approach achieves optimal average performance, it does not adequately address situations where the cost distribution exhibits heavy tails or rare but severe events. In such contexts, designing a controller that optimizes only the average performance may still allow rare but undesirable catastrophic events \cite{wang_risk-averse_2022}.

To enhance robustness, we instead seek to minimize a risk measure applied to $J(V)$. Risk measures provide a way to penalize low-probability but catastrophic scenarios. In this work, we focus on \textit{coherent risk measures}, a broad and well-established class \cite{shapiro_lectures_2009} that includes, in particular, the Conditional Value-at-Risk (CVaR). Minimizing such a measure yields controllers that reduce exposure to extreme unfavorable events.

\begin{definition} \label{def:risk}
A finite risk measure $\rho$ is a mapping from a space $\mathcal{Z} = L^p(\Omega, \mathbb{R})$, with $p \in [1, +\infty)$, to $\mathbb{R}$. A finite risk measure is said to be \textit{coherent} if it verifies the following properties:
\vspace{0em}
\begin{itemize}
\item  Convexity, $ \forall \lambda \in [0,1], \ \rho( \lambda Z + (1-\lambda) Z' ) \leq  \lambda \rho( Z ) + (1-\lambda) \rho( Z' ) $
\item Monotonicity: $ \quad Z \leq Z' \Rightarrow  \rho(Z) \leq \rho(Z')$,
\item Translation equivariance: $ \forall a \in \mathbb{R}, \ \rho(Z + a) = \rho(Z) + a $,
\item Positive homogeneity:  $\forall t \in \mathbb{R}^+, \quad \rho(tZ) = t\rho(Z)$.
\end{itemize}
\end{definition}
In what follows, we let $\rho : L^2(\Omega, \mathbb{R}) \rightarrow \mathbb{R}$ denote a fixed coherent risk measure. Our objective is to approximately solve the optimization problem
\begin{equation}\label{eq:risk_averse_problem_original}
\min_{V \in \mathcal{U}} \ \rho\bigl(J(V)\bigr),
\end{equation}
where $\mathcal{U}$ denotes the admissible control space.

\subsection{Method overview}

To approximate a solution to problem \eqref{eq:risk_averse_problem_original}, we combine the optimal control framework for the PDE–SDE interconnected system \cite{velho_optimal_2025} with a gradient-based risk-averse control method for SDEs \cite{velho_gradient_2025}, yielding the following three-step procedure.
\begin{enumerate}
\item Reformulate the coupled PDE–-SDE system as a single stochastic partial differential equation (SPDE),
\item Approximate the resulting SPDE by a finite-dimensional SDE system,
\item Apply a gradient--based risk--averse control method to the finite--dimensional SDE.
\end{enumerate}
The control obtained from Step 3 is then used as an approximate solution to the original risk--averse problem, in the same spirit as the approximation used in the mean--cost minimization setting.



\section{Approximation of the system by a SDE}\label{CC}

In this section, we first reformulate system \eqref{eq:systeme_original_coupled} as a SPDE and then apply a Galerkin projection to obtain a finite--dimensional dynamics.

\subsection{Recasting into a SPDE}

Writing the interconnected system in the SPDE framework is not straightforward, as the boundary conditions cannot be expressed as bounded operators. We therefore perform a change of variable in the PDE state as in \cite{katz_finite-dimensional_2021}, that inserts the boundary conditions inside the domain, at the cost of requiring additional regularity on the control input $V$.

Let us consider the new state variable 
\begin{equation}\label{eq:change_of_var_PDE_parabolic}
z(x,t) \triangleq u(t,x) - \theta(x) V(t) - \psi(x) M X_t
\end{equation}
where $\theta$ and $\psi$ verify the PDEs

\begin{equation}\label{eq:elliptic_PDE_rho}
 \left\{ \begin{array}{l} 
    - \Delta \theta(x)  -  c \theta(x) + \mu \theta(x) = 0 \\
\theta'(1) + \beta_1 \theta(1) = 0 \\
\theta'(0) - \beta_0 \theta(0) = 1 \\
\end{array}  \right.
\end{equation}
and
\begin{equation}\label{eq:elliptic_PDE_psi}
 \left\{ \begin{array}{l} 
    - \Delta \psi(x)  -  c \psi(x) + \mu \psi(x) = 0 \\
\psi'(1) + \beta_1 \psi(1) = 1 \\
\psi'(0) - \beta_0 \psi(0) = 0 
\end{array}  \right.
\end{equation}
For a certain value of $\mu$ sufficiently large, we can guarantee existence and uniqueness of a solution of \eqref{eq:elliptic_PDE_rho} and \eqref{eq:elliptic_PDE_psi} through the Lax-Milgram theorem \cite[Corollary 5.8]{brezis_functional_2011}. In practice, we want to choose $\mu$ relatively close to $c$, as the change of variable \eqref{eq:dynamic_change_variable_control_U_parabolic_SDE} in the control is proportional to a term $e^{\mu t}$.

Using equations \eqref{eq:elliptic_PDE_rho} and \eqref{eq:elliptic_PDE_psi}, we can compute the SPDE satisfied by $z(t,x)$:
\begin{equation}\label{eq:PDE_parabolic_SDE_z_change_var}
\left \{ \begin{array}{l}
dz(t,x) - \Delta z(t,x) dt = \mu \theta(x) V(t) dt - \theta(x) dV(t) \\
\hspace{8em} + \mu \psi(x) M X_t dt - \psi(x) M dX_t  \\
\partial_x z(t,1) + \beta_1 z(t,1) = 0 \\
\partial_x z(t,0) - \beta_0 z(t,0) =  0  \\
z(0,x) = u_0(x) - \theta(x) V(0) - \psi(x) X_0.
\end{array}
\right. 
\end{equation}
Note that we do not have directly the control inside the domain, but rather its differential $dV(t)$. Indeed, in order to actuate the PDE inside the domain, we require more regularity on the control input. 
From now on, we consider inputs $V(\cdot)$ that are $H^1(0,T)$, and we consider that our new control variable is $\Vdom$, the differential of $V$. We can then retrieve $V$ through
\begin{equation}\label{eq:dynamic_change_variable_control_U_parabolic_SDE}
\left\{ \begin{array}{l}
    dV(t) = \mu V(t) dt + \Vdom(t) dt \\
    U(0) = V_0,
\end{array} \right.
\end{equation}
with $\Vdom$ the new control variable and $V_0$ a parameter set to verify the initial compatibility condition. Let us now write the dynamic verified by the two states $(X,z)$ in the general SPDE framework. In what follows, we denote by $Z_t \triangleq \left( \begin{array}{c} X_t \\ Y_t \\ z_t \end{array} \right) \in \mathcal{H}$, the augmented state, where $\mathcal{H} \triangleq \mathbb{R}^d \times \mathbb{R} \times H^1$ the augmented state space.

\begin{proposition}
The augmented state $Z_t$ verifies the following SPDE
\small
\begin{equation}\label{eq:dynamic_SPDE_augmented_state_parabolic_SDE}
\left\{ \begin{array}{l}
dZ_t = \bigl[ ( \Tilde{\Delta} + \Tilde{A} + \Tilde{r}(t) ) Z_t + \Tilde{B} \Vdom(t) \bigr]dt + \bigl[ \Tilde{C}(t) Z_t  + \Tilde{\sigma}(t) \Bigr] dW_t \\
Z_0  =  \left( X_0, V_0,  u_0 - \theta V_0 - \psi M X_0 \right)^\top ,
\end{array} \right.
\end{equation}
\normalsize
with
\small
\begin{align*}
\Tilde{\Delta} & \triangleq \left( \begin{array}{ccc} 0 & 0 & 0 \\ 0 & 0 & 0 \\ 0 & 0 & \Delta \end{array} \right),    \text{Dom}(\Tilde{\Delta}) \triangleq \left \{  \left(  \begin{array}{c} X \\ Y \\ z  \end{array} \right) \in \mathcal{H}, z \in  \text{Dom}(\Delta)   \right \} \\ & \ \\
\Tilde{B} & \triangleq  \left( \begin{array}{c} 0 \\ 1 \\ - \theta \end{array} \right),   \Tilde{r}(t) \triangleq  \left( \begin{array}{c} r(t) \\ 0 \\ - \psi M r(t)\end{array} \right),   \Tilde{\sigma}(t) \triangleq  \left( \begin{array}{c} \sigma(t) \\ 0 \\ - \psi M \sigma(t)\end{array} \right)
\end{align*}
\begin{align*}
\Tilde{A} & \triangleq \left( \begin{array}{ccc} A - B \psi(0) M & -B\theta(0) & B \gamma_0^* \\ 0 & \mu & 0 \\ -\psi M \bigl( A - B \psi(0) M + \mu \bigr) & \psi M B \theta(0)  & - \psi M B \gamma_0^*  \end{array} \right) \\
& \ \\ 
\Tilde{C} & \triangleq \left( \begin{array}{ccc} C - D \psi(0) M & -D\theta(0) & D \gamma_0^*  \\ 0 & 0 & 0 \\ -\psi M \bigl( C - D \psi(0) M \bigr) & \psi M D \theta(0)  & - \psi M D \gamma_0^* \end{array} \right)
\end{align*}
\normalsize
where $\gamma_0 \in H^1$ is the Riesz representation of the continuous linear form
\begin{equation}\label{eq:linear_form_gamma_0_parabolic_H1_Riescz_trace_representative}
\left\{ \begin{array}{c}
\gamma_0^* : H^1 \rightarrow \mathbb{R} \\
\hspace{2em} : v \mapsto v(0).
\end{array} \right.
\end{equation}

\end{proposition}

By using the $H^1(0,T;\mathbb{R})$ norm in $\mathcal{U}$, the cost \eqref{eq:original_LQ_cost_constant} can be expressed as
\small
\begin{equation}\label{eq:original_LQ_cost_constant_SPDE_notation}
\begin{split}
J(U) \triangleq \int_0^T < \Tilde{Q} Z_t , Z_t>_{\mathcal{H}} + r U(t)^2 dt + < \Tilde{G} Z_T , Z_T>_{\mathcal{H}} 
\end{split}
\end{equation}
\normalsize
where $\Tilde{Q}$ and $\Tilde{G}$ are positive self-adjoint operators defined by
$
\Tilde{Q} \triangleq \left( \begin{array}{ccc} Q & 0 & 0 \\ 0 & r & 0 \\ 0 & 0 & 0 \end{array} \right) , \quad \Tilde{G} \triangleq \left( \begin{array}{ccc} G & 0 & 0 \\ 0 & 0 & 0 \\ 0 & 0 & 0 \end{array} \right).
$
Note that the cost is now expressed in terms of $U$ only.

\subsection{Finite dimensional approximation}

We consider the Laplacian eigenfunctions $\phi_n$ with eigenvalues $\lambda_n$, 
verifying
\begin{equation}\label{eq:stationary_PDE_eigenvector_laplacian}
\left\{ \begin{array}{l} - \Delta \phi_n(x)  = \lambda_n \phi_n \\ 
\partial_x \phi_n (0) - \beta_0 \phi_n(0) = 0, \\
\partial_x \phi_n(1) + \beta_1 \phi_n(1) = 0 ,
\end{array}  \right.
\end{equation}
and we arrange these eigenvalues in increasing order, i.e., $0 = \lambda_0 < \lambda_1 < \dots$.
In our case, we can explicitly solve the associated Sturm--Liouville problem. Note that the family $\phi_n$ forms a orthonormal basis in $H^1(0,1)$ \cite{brezis_functional_2011}.
In what follows, we denote the discretization of $H^1$ with $H^1_N \triangleq \text{Span}(\phi_0, \dots, \phi_N) $, and we accordingly denote by $\mathcal{H}_N \triangleq \mathbb{R}^d \times \mathbb{R} \times H^1_N$ the discretization of $\mathcal{H}$.
We denote by $P_N$ the orthogonal projection of $\mathcal{H}$ onto $\mathcal{H}_N$, which is given by projecting onto the first $N$ eigenvectors $\phi_0, \dots, \phi_N$. The explicit expression of $P_N$ is
\begin{equation}\label{eq:explicit_expression_projector_P_N}
    P_N \left( \begin{array}{c} X \\  Y \\ z \end{array} \right) = \left( \begin{array}{c} X \\ Y \\  \sum_{n=0}^N < z , \phi_n >_{H^1} \phi_n \end{array} \right) .
\end{equation}
The adjoint operator $P_N^*$ is the injection of $\mathcal{H}_N$ back into $\mathcal{H}$. 
We can now project the dynamic of $Z_t$ in $\mathcal{H}_N$.
We define the finite dimensional matrices $\Tilde{\Delta}_N = P_N \Tilde{\Delta} P_N^* $, $\Tilde{A}_N = P_N \Tilde{A} P_N^* $ $\Tilde{C}_N = P_N \Tilde{C} P_N^* $, $\Tilde{B}_N = P_N \Tilde{B} $, $\Tilde{r}_N = P_N \Tilde{r} $ and $\Tilde{\sigma}_N = P_N \Tilde{\sigma} $.
The approximated state $Z_t^N \in \mathcal{H}_N$ is the solution to the linear SDE
\begin{equation}\label{eq:formal_SPDE_discretized_linear_system_LQ_finite_dim}
\left \{ \begin{array}{l}
dZ_t^N = \bigl[ ( \Tilde{\Delta}_N + \Tilde{A}_N ) Z_t^N + \Tilde{B}_N U(t) + \Tilde{r}_N(t) \bigr] dt \\
\hspace{3em} + \bigl[ \Tilde{C}_N Z_t^N + \Tilde{\sigma}_N(t) \bigr] dW_t,\\
Z^N_0 =  P_N Z_0 .
\end{array} \right.
\end{equation}
The different matrices can be written explicitly along the basis of vectors $(\phi_0, \dots, \phi_N)$ 
In particular, since we chose to project onto spans of eigenfunctions of $\Delta$, we have that $\Tilde{\Delta}_N \triangleq \text{diag}( 0_{d+1}, \lambda_0, \dots \lambda_N) $.

The cost functional given in \eqref{eq:original_LQ_cost_constant_SPDE_notation} is in turn approximated by 
$$
J_N(U) \triangleq \int_0^T \bigl[ (Z_t^N)^\top \Tilde{Q}_N Z_t^N + r U(t)^2 \bigr] dt +  (Z_T^N)^\top \Tilde{G}_N Z_T^N 
$$
with $\Tilde{Q}_N \triangleq P_N \Tilde{Q} P_N^* $ and $\Tilde{G}_N \triangleq P_N \Tilde{G} P_N^* $.

\subsection{Some convergence results}

In what follows, we denote by $S(\cdot)$ and $S_N(\cdot)$ the strongly continuous semi-group generated respectively by $\Tilde{\Delta}$ and $\Tilde{\Delta}_N$.
It is a known result that the spectral approximation of the Laplacian $\Tilde{\Delta}_N$ 
verifies the Trotter-Kato theorem \cite[Thm 4.3]{morris_controller_2020}. Consequently, we have that 
\begin{equation}\label{eq:Trotter_Kato_consequence_convergence_S}
\forall Z \in \mathcal{H},~ \lim_{N \rightarrow + \infty} \sup_{0 \leq t \leq T} \Vert P_N^* S_N(t) P_N Z - S(t)Z \Vert_\mathcal{H} = 0,
\end{equation}
Moreover, we can verify that all the other bounded operators converge strongly pointwise, that is
\begin{equation}\label{eq:convergence_bounded_operators}
\begin{split}
\forall Z \in \mathcal{H}, \quad & \Vert \Tilde{A} Z - P_N^* \Tilde{A}_N P_N Z \Vert_\mathcal{H} \rightarrow 0,  \\
& \Vert \Tilde{C} Z - P_N^* \Tilde{C}_N P_N Z \Vert_\mathcal{H} \rightarrow 0, \\
&  \Vert \Tilde{B} - P_N^* \Tilde{B}_N \Vert_\mathcal{H} \rightarrow 0.
\end{split}
\end{equation}
These convergences also hold for the adjoints of the operators $\Tilde{A}^*$, $\Tilde{B}^*$ and $\Tilde{C}^*$. Since the approximate dynamics converge to those of the original PDE–SDE system, the corresponding cost functionals also converge. By continuity of the risk measure, any minimizer of \(U \mapsto \rho(J_N(U))\) is therefore a consistent approximation of a solution to \eqref{eq:risk_averse_problem_original}.


\section{Gradient Method for Finite Dimensional Risk Averse Control}\label{CS}

In this section, we adapt the method developed in \cite{velho_gradient_2025} to compute an approximate minimizer of the mapping $U \in \mathcal{U} \mapsto \rho(J_N(U))$. We begin by specifying a practical choice of the control space $\mathcal{U}$. We then reformulate the original nonlinear minimization problem as an equivalent linear–convex min–max problem. Finally, we describe a gradient--based algorithm for computing a saddle point of this min–max formulation.

\subsection{Choice of control space}\label{CHOICE_CONTROL}

The space of admissible controls is $(v,K) \in L^2([0,T] , \mathcal{C}) \times \mathcal{K}$, with $\mathcal{C} \subset \mathbb{R}$ and $\mathcal{K} \subset \mathbb{R}^{1 \times (d+1+N)}$  two compact convex sets.
More specifically, to each control parameter $ (v, K) \in \mathcal{U}$ we associate a feedback control $U(\cdot)$ in the following manner:
\begin{equation}\label{eq:markovian_control_feedback_formula_definition_risk_averse_finDim}
    U_{v, K}(t) = v(t) + K Z_t^N,
\end{equation}
where $Z_t^N$ follows the discretized dynamic \eqref{eq:formal_SPDE_discretized_linear_system_LQ_finite_dim}.
This control structure, commonly used in stochastic control \cite{ lew_sample_2024, lew_risk-averse_2024}, combines a deterministic open--loop term $v(\cdot)$ with a state-feedback gain $K$. The open--loop part steers the nominal system behavior, while the feedback term adapts to stochastic fluctuations through real--time state measurements. This avoids the numerical difficulties of introducing randomness directly into $v(\cdot)$ \cite{kushner_numerical_2001}, since noise is accessed implicitly via the observed state.
Our objective is to determine the parameters $(v,K)$ that minimize the quantity $\rho\bigl(J_N(U_{v,K})\bigr)$.

\subsection{Recasting into a min--max problem}

The functional $U \mapsto \rho(J_N(U))$ is in general nonlinear and costly to evaluate directly. To handle this, we rely on a fundamental representation of coherent risk measures, which expresses them as the supremum of linear expectations \cite{shapiro_lectures_2009}.

\begin{proposition}
Let $\rho$ be a coherent risk measure, then
\begin{equation}
\forall Z \in L^2(\Omega), \quad     \rho(Z) = \max_{\zeta \in \partial \rho(0)} \ \espE \bigl[  \zeta Z  \bigr]
\end{equation}
where $\partial \rho(0)$ denotes the subgradient of $\rho$ at 0.
\end{proposition}

Using this representation, we can recast the minimization problem into
$$
\min_{U \in \mathcal{U}} \rho( J_N(U) ) = \min_{U \in \mathcal{U}} \ \max_{\zeta \in \partial \rho(0)} \espE[  \zeta  \, J_N(U) ].
$$
In the remainder of this section, we present a gradient-based procedure to compute a saddle point of the function $(\zeta, U) \mapsto \mathbb{E}\bigl[\zeta \, J_N(U)\bigr]$.

\subsection{A gradient ascent-descent algorithm}

Our projected gradient descent-ascent algorithm updates at each iteration $n$ our estimates $\zeta_n$ and $(v_n , K_{n})$ with the following rule: 
\small
\begin{equation}\label{eq:update_rule_modified_algo_risk_averse_finDim}
\begin{split}
    &v_{n+1} = \mathcal{P}_{ L^2([0,T], \mathcal{C}) } \Bigl(  v_n - \eta \espE \bigl[ \zeta_n (1 - \gamma \zeta_n) \nabla_{v} J( U_{v_n , K_{n}}  ) \bigr] \Bigr) , \\
    & K_{n+1} = \mathcal{P}_{\mathcal{K}} \Bigl( K_{n} - \eta \espE \bigl[ \zeta_n (1 - \gamma \zeta_n) \nabla_{K} J( U_{v_n , K_{n}}  ) \bigr] \Bigr), \\
    & \zeta_{n+1} = \mathcal{P}_{\partial \rho (0)} \Bigl( \zeta_n + \beta (1 - 2\gamma \zeta_n) J( U_{v_n , K_{n}} ) \Bigr).
\end{split}
\end{equation}
\normalsize
where $\gamma > 0$ is a regularization parameter, taken small, $ \mathcal{P}_{K}$ denotes the projection on a closed convex set $K$, and $\eta$ and $\beta$ are two positive constants that represent the step size of the gradient. To increase convergence speed, one could make these steps vary from iteration to iteration based on the change of direction of the gradient \cite{kingma_adam_2017}. 


The updates  \( ( v_{n+1}, K_{n+1}) \) must be projected onto the closed convex control set, since gradient steps need not remain feasible. The same projection step applies to \( \zeta_{n+1} \). In practice, these projections are computationally inexpensive when \( \mathcal{U} \) is chosen with simple constraints.

The practical implementation of \eqref{eq:update_rule_modified_algo_risk_averse_finDim} requires computing the gradient of the cost $J_N$ with respect to $v$ and $K$. Below, we summarize the expressions of this gradient, referring the reader to \cite{velho_gradient_2025} for their proof. 

\begin{proposition}
Fix a control $U_{v, K} \in \mathcal{U}$. The Riesz representation of the gradient of $J$ at $U_{v, K}$ is 
\begin{align*}
&\nabla_v J_N(U_{v, K})(t) = \Bigl( I_U(t) \ \phi_{U}^{-1}(t) \ \Tilde{B}_N \Bigr)^\top + 2 r U_{v,K}(t) \\
&\nabla_K J_N(U_{v, K}) = \int_0^T \Bigl[ \Bigl( I_U(t) \ \phi_{U}^{-1}(t) \ \Tilde{B}_N \Bigr)^\top Z^N_{U}(t)^\top  \\ 
& \hspace{12em} +  2 r U_{v,K}(t)  Z^N_{U}(t)^\top \Bigr] dt
\end{align*}
where $Z^N_{U}(t)$ follows the dynamic \eqref{eq:formal_SPDE_discretized_linear_system_LQ_finite_dim} associated to the control $U_{v, K}$, and
\small
\begin{equation*}
\begin{split}
I_U(t) \triangleq \int_t^T \Bigl( Z^N_{U}(s)^\top \Tilde{Q}_N  +  2 r U_{v,K}(s) K \Bigr)  \phi_{U}(s) ds +  Z^N_{U}(T)^\top \Tilde{G}_N
\end{split}
\end{equation*}
\normalsize
Here, $\phi_U(.)$ is the fundamental matrix solution to the closed loop SDE, i.e., the unique solution to
\begin{equation}\label{eq:linearized_sde_phi_risk_averse_finDim}
    \left\{
    \begin{array}{l}
    d\phi_U(t) = \Bigl [ \Tilde{\Delta}_N + \Tilde{A}_N +  \Tilde{B}_N K \Bigr ] \phi_U(t) dt + \Tilde{C}_N \phi_U(t) dW_t,\\
    \phi_U(0) = Id_{d+N+1} .
    \end{array}
    \right.
\end{equation}
\end{proposition}

Thus, the primary challenge of implementing this method lies in simulating the trajectories of $X_u$ and $\phi_u$ and computing expectations of random variables that involve terms from those trajectories.



\section{Numerical results}\label{NR}



We illustrate our method on an academic example. We consider the unstable interconnected system
\begin{equation}\label{eq:systeme_numerics_considered}
\left\{ \begin{array}{l} \partial_t u(t,x) = \Delta u(t,x)  +  0.2 u(t,x) \\ 
dX_t = \bigl( A X_t + B u(t,0) \bigr) dt + \bigl( C X_t + \sigma(t) \bigr) dW_t \\
\partial_x u(t,1) = V(t) \\
\partial_x u(t,0) = 0
\end{array}  \right.
\end{equation}
with $A = \left( \begin{array}{cc} 0.6 &  0.4  \\ 0 &  0.4 \end{array} \right)$, $B = \left( \begin{array}{c} 0  \\ 1 \end{array} \right)$, $\sigma(t) = \left( \begin{array}{c} 0.05  \\ 0.05 \end{array} \right)$, $C = \text{diag}(0.1, 0.1)$ and $T = 4$.

We minimize the quadratic cost
\begin{equation}
\begin{split}
J(V) & = c_q \int_0^T \Vert X(t) \Vert^2 dt +  r \int_0^T \bigl(  V(t)^2 + V'(t)^2 \bigr) dt . 
\end{split}
\end{equation}
with $c_q = 1$ and $r = 3$.
For the Galerkin approximation, we take \(N=3\), since the eigenvalues \(\lambda_n = n^2 \pi^2\) grow rapidly and the higher modes decay as \(e^{(c-\lambda_n)t}\), becoming negligible for \(n>3\). 
For the simulation of the system, we use a Galerkin method (on a different basis of functions) for the PDE, and a Euler--Maruyama scheme for the SDE.

As a risk measure, we use the Conditional Value--at--Risk \(CVaR_\alpha\), a widely employed coherent risk measure \cite{shapiro_lectures_2009}. Intuitively, \(CVaR_\alpha(Z)\) corresponds to the mean of the worst \(\alpha\)--fraction of outcomes. We choose \(\alpha = 0.1\), so that \( CVaR_{0.1}(J(U)) \) focuses on the 10$\%$ highest--cost scenarios. Its subdifferential at zero is
$$
\partial CVaR_\alpha(0) = \bigl\{\zeta \in L^2(\Omega): \mathbb{E}[\zeta]=1,\ \zeta \in [0,\alpha^{-1}] \text{ a.s.}\bigr\},
$$
for which projection is numerically inexpensive.

We initialize the control using the feedback law that minimizes the expected cost \(\mathbb{E}[J_N(U)]\), obtained from standard stochastic LQ theory \cite{yong_stochastic_1999} and compute the gain numerically through a Kleinman-Newton method.


\begin{figure}[ht!]
    \centering
    \includegraphics[width=0.95\linewidth]{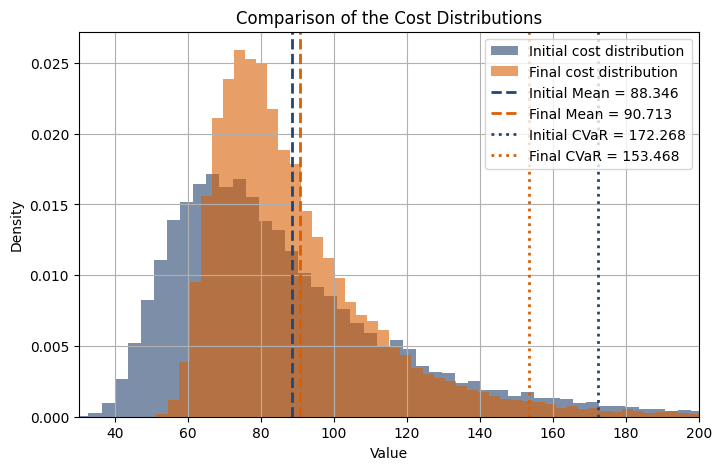}
    \caption{Comparison of initial and final distribution of the cost on $10000$ realizations of the trajectories}
    \label{fig:output_final_distribs_jolies}
\end{figure}

We then compare the cost distribution obtained with the mean--optimal control (the initial guess) to the distribution resulting from our algorithm after 1000 iterations.
Figure \ref{fig:output_final_distribs_jolies} shows that, while the average performance decreases slightly, the upper tail of the cost distribution is significantly reduced. This leads to more reliable behavior, reflected in a markedly lower risk measure.

\begin{table}[ht!]
\centering
\begin{tabular}{c c c c}
\hline
\textbf{Quantile} & \textbf{Cost 1} & \textbf{Cost 2} & \textbf{Difference (Cost1 - Cost2)} \\
\hline
0.20 & 60.90 & 71.35 & -10.45 \\
0.40 & 72.94 & 79.31 & -6.37 \\
0.60 & 86.58 & 88.59 & -2.01 \\
0.70 & 96.03 & 95.17 & 0.86 \\
0.80 & 109.27 & 104.72 & 4.55 \\
0.90 & 132.04 & 121.37 & 10.67 \\
0.95 & 157.51 & 140.74 & 16.76 \\
0.99 & 222.28 & 193.27 & 29.01 \\
\hline
\end{tabular}
\caption{Comparison of quantiles of the cost distributions.}
\label{tab:quantile_comparison}
\end{table}

The quantile values provide a more quantitative comparison between the two controls. In particular, for the worst $10\%$ of scenarios, Table \ref{tab:quantile_comparison} shows that the risk--averse controller yields substantially lower costs, confirming its improved robustness to rare but severe events.


\section{Conclusion}\label{CL}

In this paper, we proposed a method for risk--averse optimal control of an interconnected system consisting of a heat PDE coupled with a linear SDE. The approach extends control strategies that minimize the expected value of the cost by instead optimizing a coherent risk measure, thereby improving robustness to rare but high-impact events. Numerical results demonstrate that the method significantly reduces the tail of the cost distribution while maintaining reasonable average performance.

Future research will focus on establishing theoretical guarantees for the complete procedure by combining approximation results for the PDE-–SDE reduction with convergence results for the gradient-based risk--averse optimization. This would provide convergence bounds for the overall method. Additionally, we aim to extend the framework to systems with nonlinear stochastic dynamics, further broadening its applicability.

\bibliographystyle{ieeetr}
\bibliography{references}

\begin{thebibliography}{10}

\bibitem{tang_state_2011}
S.~Tang and C.~Xie, ``State and output feedback boundary control for a coupled {PDE}–{ODE} system,'' {\em Systems \& Control Letters}, vol.~60, pp.~540--545, Aug. 2011.

\bibitem{antonio_susto_control_2010}
G.~Antonio~Susto and M.~Krstic, ``Control of {PDE}–{ODE} cascades with {Neumann} interconnections,'' {\em Journal of the Franklin Institute}, vol.~347, pp.~284--314, Feb. 2010.

\bibitem{mohammadi_optimal_2015}
L.~Mohammadi, I.~Aksikas, S.~Dubljevic, and J.~F. Forbes, ``Optimal boundary control of coupled parabolic {PDE}–{ODE} systems using infinite-dimensional representation,'' {\em Journal of Process Control}, vol.~33, pp.~102--111, Sept. 2015.

\bibitem{auriol_late-lumping_2019}
J.~Auriol, K.~A. Morris, and F.~Di~Meglio, ``Late-lumping backstepping control of partial differential equations,'' {\em Automatica}, vol.~100, pp.~247--259, Feb. 2019.

\bibitem{vazquez_backstepping_2026}
R.~Vazquez, J.~Auriol, F.~Bribiesca-Argomedo, and M.~Krstic, ``Backstepping for partial differential equations: {A} survey,'' {\em Automatica}, vol.~183, p.~112572, Jan. 2026.

\bibitem{touzi_optimal_2013}
N.~Touzi, {\em Optimal {Stochastic} {Control}, {Stochastic} {Target} {Problems}, and {Backward} {SDE}}, vol.~29 of {\em Fields {Institute} {Monographs}}.
\newblock New York, NY: Springer, 2013.

\bibitem{lew_risk-averse_2024}
T.~Lew, R.~Bonalli, and M.~Pavone, ``Risk-{Averse} {Trajectory} {Optimization} via {Sample} {Average} {Approximation},'' {\em IEEE Robotics and Automation Letters}, vol.~9, pp.~1500--1507, Feb. 2024.
\newblock Conference Name: IEEE Robotics and Automation Letters.

\bibitem{talebi_review_2016}
B.~Talebi, P.~A. Mirzaei, A.~Bastani, and F.~Haghighat, ``A {Review} of {District} {Heating} {Systems}: {Modeling} and {Optimization},'' {\em Frontiers in Built Environment}, vol.~2, Oct. 2016.

\bibitem{velho_optimal_2025}
G.~Velho, J.~Auriol, I.~Boussaada, and R.~Bonalli, ``Optimal {Control} of an {Interconnected} {SDE} -{Parabolic} {PDE} {System}.'' CDC 2025 - Submit status: Accepted, Apr. 2025.

\bibitem{wang_risk-averse_2022}
Y.~Wang and M.~P. Chapman, ``Risk-averse autonomous systems: {A} brief history and recent developments from the perspective of optimal control,'' {\em Artificial Intelligence}, vol.~311, p.~103743, Oct. 2022.

\bibitem{shapiro_lectures_2009}
A.~Shapiro, D.~Dentcheva, and A.~Ruszczyński, {\em Lectures on {Stochastic} {Programming}: {Modeling} and {Theory}}.
\newblock Society for Industrial and Applied Mathematics, Jan. 2009.

\bibitem{chapman_risk-sensitive_2022}
M.~P. Chapman, R.~Bonalli, K.~M. Smith, I.~Yang, M.~Pavone, and C.~J. Tomlin, ``Risk-{Sensitive} {Safety} {Analysis} {Using} {Conditional} {Value}-at-{Risk},'' {\em IEEE Transactions on Automatic Control}, vol.~67, pp.~6521--6536, Dec. 2022.
\newblock Conference Name: IEEE Transactions on Automatic Control.

\bibitem{bonalli_first-order_2023}
R.~Bonalli and B.~Bonnet, ``First-{Order} {Pontryagin} {Maximum} {Principle} for {Risk}-{Averse} {Stochastic} {Optimal} {Control} {Problems},'' May 2023.
\newblock arXiv:2204.03036 [math] version: 6.

\bibitem{miller_optimal_2017}
C.~W. Miller and I.~Yang, ``Optimal {Control} of {Conditional} {Value}-at-{Risk} in {Continuous} {Time},'' {\em SIAM Journal on Control and Optimization}, vol.~55, pp.~856--884, Jan. 2017.

\bibitem{velho_gradient_2025}
G.~Velho, J.~Auriol, and R.~Bonalli, ``A gradient descent-ascent method for continuous-time risk-averse optimal control,'' {\em Systems \& Control Letters}, vol.~204, p.~106186, Oct. 2025.

\bibitem{brezis_functional_2011}
H.~Brezis, {\em Functional {Analysis}, {Sobolev} {Spaces} and {Partial} {Differential} {Equations}}.
\newblock New York, NY: Springer, 2011.

\bibitem{katz_finite-dimensional_2021}
R.~Katz and E.~Fridman, ``Finite-dimensional control of the heat equation: {Dirichlet} actuation and point measurement,'' {\em European Journal of Control}, vol.~62, pp.~158--164, Nov. 2021.

\bibitem{morris_controller_2020}
K.~A. Morris, {\em Controller {Design} for {Distributed} {Parameter} {Systems}}.
\newblock Communications and {Control} {Engineering}, Cham: Springer International Publishing, 2020.

\bibitem{lew_sample_2024}
T.~Lew, R.~Bonalli, and M.~Pavone, ``Sample {Average} {Approximation} for {Stochastic} {Programming} with {Equality} {Constraints},'' {\em SIAM Journal on Optimization}, vol.~34, pp.~3506--3533, Dec. 2024.
\newblock Publisher: Society for Industrial and Applied Mathematics.

\bibitem{kushner_numerical_2001}
H.~J. Kushner and P.~Dupuis, {\em Numerical {Methods} for {Stochastic} {Control} {Problems} in {Continuous} {Time}}, vol.~24 of {\em Stochastic {Modelling} and {Applied} {Probability}}.
\newblock New York, NY: Springer, 2001.

\bibitem{kingma_adam_2017}
D.~P. Kingma and J.~Ba, ``Adam: {A} {Method} for {Stochastic} {Optimization},'' Jan. 2017.
\newblock arXiv:1412.6980 [cs].

\bibitem{yong_stochastic_1999}
J.~Yong and X.~Y. Zhou, {\em Stochastic controls: {Hamiltonian} systems and {HJB} equations}.
\newblock No.~43 in Applications of mathematics, New York: Springer, 1999.

\end{thebibliography}

\end{document}